\documentclass[12pt]{article}
\input{amssym.def}
\input{amssym}
\usepackage{eufrak}

\input{amssym.def}
\input{amssym}

\def\theequation{\thesection.\arabic{equation}}
\makeatletter \@addtoreset{equation}{section} \makeatother
\def\nn{\nonumber}
\def\be{\begin{equation}}
\def\ee{\end{equation}}
\def\ba{\begin{eqnarray}}
\def\ea{\end{eqnarray}}
\def\lb{\label}

\def\tbl#1#2{{\ifmmode
\Bigl\{\!\!\begin{array}{c}
\scriptstyle #1\\[-2pt]
\raisebox{2pt}{$\scriptstyle #2$} \end{array}\!\!\Bigr\} \else
\raisebox{2pt}{\scriptsize$\left\{\!\!\!\begin{array}{c}
#1\\[-1pt] #2 \end{array}\!\!\!\right\}$}\fi}}

\newtheorem{theorem}{Theorem}

\newtheorem{conjecture}[theorem]{Conjecture}

\newtheorem{proposition}[theorem]{Proposition}

\begin{document}

\def\h{{\hbar}}
\def\O{\cal O}
\def\Cas{{\rm Cas}}
\def\qq{q^{-1}}
\def\gq{\gggg_{q}}
\def\uq{U_q(sl(2))}
\def\Uq{U_q(sl(n))}
\def\Uqq{U_q(\gggg)}
\def\C{{\Bbb C}}
\def\R{{\Bbb R}}
\def\K{{\Bbb K}}
\def\De{\Delta}
\def\de{\delta}
\def\ve{\varepsilon}
\def\ep{\epsilon}
\def\al{{\alpha}}
\def\ot{\otimes}
\def\om{\omega}
\def\Om{{\Omega}}
\def\End{{\rm End\, }}
\def\Tr{{\rm Tr}}
\def\Sym{{\rm Sym\, }}
\def\Gr{{\rm Gr}}
\def\dim{{\rm dim}}
\def\sdim{{\rm sdim}}
\def\vv{V^{\ot 2}}
\def\vl{V_{\la}}
\def\la{{\lambda}}
\def\span{\rm span}
\def\lhq{{\cal L}(R,\h)}
\def\lqh{{\cal L}(R,\h)}
\def\lq{{\cal L}(R)}
\def\lh{{\cal L}(\h)}
\def\Omn{\C_q[{\cal O}_{\mu, \nu}]}
\def\M{{\cal M}}
\def\L{{\cal L}}
\def\LL{{\cal L}^{\ot 2}}
\def\sla{s_{\la}}
\def\schf#1{s_{\raisebox{-0.08mm}{$_{#1}$}}}
\def\E{\cal E}
\def\hE{\widehat{\cal E}}
\def\oe{\overline{e}}
\def\he{\hat{e}}
\def\oM{\overline{M}}
\def\hL{\widehat{L}}
\def\hl{\widehat{l}}
\def\ha{\widehat{a}}
\def\hc{\widehat{c}}
\def\hm{\widehat{\mu}}
\def\hn{\widehat{\nu}}
\def\hmu{\widehat{\mu}}
\def\hnu{\widehat{\nu}}
\def\hd{\widehat{d}}
\def\hp{\widehat{p}}
\def\hOmn{\C_q[\widehat{{\cal O}_{\mu, \nu}}]}
\def\Mmn{M_{\mu,\nu}}
\def\hMmn{\widehat{M_{\mu,\nu}}}

\def\bea{\begin{eqnarray}}
\def\eea{\end{eqnarray}}
\def\be{\begin{equation}}
\def\ee{\end{equation}}
\def\nn{\nonumber}

\makeatletter
\renewcommand{\theequation}{{\thesection}.{\arabic{equation}}}
\@addtoreset{equation}{section} \makeatother

\title{Generic super-orbits in $gl(m|n)^*$ and their braided counterparts}
\author{
\rule{0pt}{7mm} Dimitri
Gurevich\thanks{gurevich@univ-valenciennes.fr}\\
{\small\it LAMAV, Universit\'e de Valenciennes,
59313 Valenciennes, France}\\
\rule{0pt}{7mm} Pavel Saponov\thanks{Pavel.Saponov@ihep.ru}\\
{\small\it Division of Theoretical Physics, IHEP, 142281 Protvino,
Russia} }

\maketitle

\begin{abstract}
We introduce some braided varieties --- braided orbits --- by
considering quotients of the so-called Reflection Equation Algebras associated
with Hecke symmetries (i.e. special type solutions of the quantum
Yang-Baxter equation). Such a braided variety is called regular if
there exists a  projective module on it,  which is a counterpart of
the cotangent bundle on a generic orbit ${\O}\in gl(m)^*$ in the
framework of the Serre approach. We give a criterium of regularity of
a braided orbit in terms of roots of the Cayley-Hamilton identity
valid for the generating matrix of the Reflection Equation Algebra
in question. By specializing our general construction we get
super-orbits in $gl(m|n)^*$ and a criterium of their regularity.
\end{abstract}

{\bf AMS Mathematics Subject Classification, 2010:} 81R60


\section{Introduction}

Let  $L\in gl(m,\C)^*$ be a generic matrix, that is its eigenvalues $\mu_i$, $1\leq i \leq m$,
are pairwise distinct. As is well known, its orbit
${\O} \subset gl(m,\C)^*$ under the coadjoint action of the group
$GL(m,\C)$\footnote{In what follows we omit the sign $\C$ in
our notation.} is a regular affine algebraic variety. Namely, any
such  orbit $\O$ can be defined by the following system of polynomial
equations:
\be
\Tr L^k-\sum_{i=1}^m \mu_i^k=0,\quad k=1,...,m.
\lb{sys}
\ee

The main objective of the present paper is to give an explicit
description of  generic orbits in $gl(m|n)^*$ and their braided
counterparts. These braided generic orbits constitute a
subclass of braided varieties  which are quotients of the
so-called Reflection Equation Algebra of $GL(m|n)$ type. Let us
recall its definition.

Let $V$ be a finite dimensional vector space. A {\em braiding} $R\in \End(\vv)$ is a solution of
the quantum Yang-Baxter equation
$$
R_{12} R_{23}R_{12}=R_{23}R_{12}R_{23},\quad{\rm where}\quad R_{12}=R\ot I,\,
R_{23}=I\ot R.
$$
A braiding $R$ is called  {\it a Hecke symmetry} if it satisfies the following second degree equation
$$
(qI-R)(\qq I+R)=0,
$$
where $q\in \C^{\times}$ is generic. If $q=1$ (this value is not forbidden) the corresponding braiding
is called an {\em involutive symmetry}.

Given a Hecke symmetry $R$, consider a unital
associative algebra    generated by elements $l^i_j$ subject to the system
\be
R\,L_1\,R\,L_1-L_1\,R\,L_1\, R=0,
\label{RE}
\ee
where $L_1=L\ot I$ and $L=\|l_i^{j}\|$, $1\leq i,\,j\leq \dim\,V$,
is a matrix with entries $l_i^{j}$. This algebra is denoted $\lq$ and called {\em the Reflection Equation Algebra}
(REA) corresponding to the given Hecke symmetry $R$,  the matrix $L$  is called the {\em generating} matrix
of the algebra $\lq$.

We say that an algebra $\lq$ is of $GL(m|n)$  type if the degree of the numerator (resp., denominator)
of the Hilbert-Poincar\'e series $P_-(t)$  (see section \ref{sec:2}), which is always a rational function, equals $m$ (resp., $n$).
A popular example is provided by the algebra $\lq$ corresponding to a Hecke symmetry $R=R(q)$ which is a deformation
of a super-flip $\sigma=R(1)$. The operator $\sigma$ acts in the space $\vv$ where $V=V_0\oplus V_1$ is a super-space with the even
component $V_0$ and odd component $V_1$ such that $\dim V_0=m$ and $\dim V_1=n$ (the
ordered couple  $(m|n)$ is usually called  the super-dimension\footnote{In general, we call the couple $(m|n)$ the bi-rank. Also,
note that in the notation $\Sym(gl(m|n))$ below  the symmetric algebra is understood in the sense of the
super-theory.} of the space $V$). In this case the algebra $\lq$ is a deformation of $\C[gl(m|n)^*]\cong \Sym(gl(m|n))$ and
turns into the latter algebra as  $q\to 1$. In particular,  if the space $V$ is even (i.e. $n=0$), the algebra $\lq$ is a deformation
of the algebra  $\C[gl(m)^*]$. As an example we mention the REA corresponding to the Hecke symmetry coming from the
Quantum Group $U_q(sl(m))$.

For a more detailed treatment of the $GL(m|n)$  type REA we refer the reader to \cite{GPS1, GPS2, GPS3}. Here
we reproduce one of the basic properties of this algebra: its center $Z(\lq)$ is similar to that of the enveloping algebra $U(gl(m|n))$.
In particular,  the elements
$$
p_k(L):=\Tr_R L^k, \quad k=1,2,...
$$
called the {\em power sums} in analogy with the classical case, belong to the center $Z(\lq)$.
Hereafter, $\Tr_R$ stands for the so-called  braided (quantum or $R$-) trace which can be associated with any
skew-invertible (see section \ref{sec:2}) braiding $R$. Let us emphasize that this trace is one of the main features of the braided
geometry\footnote{Another
important feature of the braided geometry is a modification of the notions of  Lie algebras,  vector fields,
differential operators. All these notions are coordinated with the initial Hecke symmetry $R$.
Thus, if $R$ is a super-flip, $\Tr_R$ turns into a (super-)trace and the corresponding  "braided Lie bracket" $[\,,\,]_R$
turns into a super-Lie one.} as defined in \cite{GS2}.

Since the elements  $\Tr_R L^k$ are central in the algebra $\lq$, it is natural
to consider the quotient $\lq/\langle I \rangle$ where $\langle I \rangle $ is the two-sided ideal generated
by the set $I\subset \lq$ consisting of the elements coming in the left hand side of  (\ref{sys}) but with $\Tr_R$
instead of the usual trace. Also, the sums  $\sum_{i=1} \mu_i^k$, $k=1,2,...,m$, must be modified in an
appropriate way. Such a quotient $\lq/\langle I \rangle$
can be treated as a braided analog of (the coordinate rings of) an algebraic variety, which is a generic
coadjoint orbit provided  all $\mu_i$ are pairwise distinct. However, proceeding in this way, we have to answer the following questions.
 \begin{enumerate}
\item
What are the "braided" analogs of the eigenvalues $\mu_i$?
\item
How many equations are there in the braided case, i.e. which number must replace  the index $m$ in (\ref{sys})?
\item
For which values of the braided eigenvalues the corresponding quotient can be treated as a {\it regular}
braided  variety (and consequently, a braided generic orbit)?
\end{enumerate}
Note that even if $R$ is a super-flip and therefore $\lq\cong \Sym(gl(m|n))$, these questions are still meaningful.

The problem of diagonalization of a super-matrix was studied in \cite{Sh}. Let us point out that we do not consider
such a diagonalization. We define the eigenvalues of a super-matrix as the roots of the Cayley-Hamilton (CH) identity
satisfied by this super-matrix. The various forms of the CH identity for super-matrices (including the one convenient for
our aims) has been given in \cite{KT}. In \cite{GPS1} a CH identity was presented for the generating
matrix\footnote{\label{ftn:4}Note, that we do not speak about the CH identity and "eigenvalues" of an
{\it arbitrary} matrix with entries from $\lq$. We are dealing with a very special  matrix $L$. In a sense,
it arises from the central element $\Tr_R L^2$ as explained in \cite{GLS}.

Also, note that the algebra $\lq$ is a particular case of the so-called Quantum Matrix Algebras (QMA) which are associated
with a couple of compatible braidings $(R,F)$ (see \cite{GPS1}). The generating matrix $L$ of a QMA satisfies a CH identity
as well, but its form differs from the classical one. Besides, in general, the coefficients of such a CH identity are not central in
the corresponding  algebra. This obstacle does not allow us to consider similar "orbits" in other QMA.}  $L$ of the REA $\lq$
associated with any skew-invertible Hecke symmetry $R$ of the $GL(m|n)$ type. The corresponding CH identity has the form
\be
\sum_{k=1}^{m+n} c_k(L) L^k=0,
\label{CH-ident}
\ee
where $c_k(L)$ are non-trivial central elements of the algebra $\lq$.

Let $\{\mu_i\}_{1\leq i\leq m+n}$ be the roots of the equation
$$
\sum_{k=1}^{m+n} c_k(L) \mu^k=0
$$
considered as elements of the algebraic extension of  the localization
$$
Z(\lq)_{loc}:=S^{-1}Z(\lq)
$$
of  $Z(\lq)$ by the set $S=\{c_{m+n}^k(L),\,k=1,2,...\}$. These roots are called
{\em quantum eigenvalues} of the matrix $L$. We assume them to be central in the algebra $\lq$.
They play the role of usual eigenvalues in the
analysis below. This is the answer to the question 1 from the above list. If $n=0$
the leading  coefficient $c_{m+n}(L)$ equals 1 and the mentioned
localization does not affect the algebra $Z(\lq)$.  However, in general, $c_{m+n}(L)$ is not a number.

Also, in the general case (i.e. if $n\not=0$) the set of all  eigenvalues splits into two subsets: even eigenvalues and
odd ones. This splitting stems from the factorization of the CH identity discovered in \cite{GPS1}.
In what follows we denote the odd eigenvalues $\nu_i$ and keep the notation $\mu_i$ for even ones.

The key point of our method is a parametrization of the power sums $\Tr_R L^k$ in terms of the quantum
eigenvalues. Let $p_k(\mu, \nu)$ be such a parametrization\footnote{In the super-case this parametrization reads
$$
p_k(\mu, \nu)=\sum_{i=1}^m \mu_i^k-\sum_{j=1}^n \nu_j^k.
$$}.
Then we define a braided variety $\Omn$
by the following system of the polynomial equations
\be
\Tr_R L^k-p_k(\bar\mu, \bar\nu)=0,\qquad k=1,2,\dots,m+n,\quad \bar\mu,\bar\nu\in\C.
\label{sys1}
\ee
Here we pass to a specialization $\mu\mapsto \bar\mu$ and $\nu\mapsto \bar\nu$ of the elements from $Z(\lq)_{loc}$
to complex numbers. In what follows we  omit the bar over the letters keeping the notations
 $\mu$ and $\nu$ for numeric values of roots.  In each case
the meaning of a symbol is clear from the context.

Thus, the braided variety $\Omn$ is a quotient algebra
$$
\Omn=\lq/\langle \Tr_R L-p_1(\mu, \nu),...,\Tr_R L^{m+n}-p_{m+n}(\mu, \nu)  \rangle.
$$
By this we give an answer to the question 2 from the above list: the number of defining equations in the system
(\ref{sys1}) equals $m+n$. By abusing the language\footnote{Note that in the classical case such a variety is not an orbit
but rather a collection of them if eigenvalues of the matrix $L$ are not pairwise distinct.}
we call this braided variety  {\em a braided orbit}.

As for the question 3, it  can be now reformulated in the following way. For which values of  $\mu=(\mu_1,...\mu_m)$
and $\nu=(\nu_1,...\nu_n)$ the quotient $\Omn$ can be considered as a regular  braided  variety?
 Our answer to this question is based on the following observation. According to the famous Serre result \cite{S}
the space of sections of a  vector bundle
on a regular affine algebraic variety $\M$ is a finitely generated projective module over the
coordinate algebra $\C[\M]$. On quantizing this algebra, it is possible to simultaneously quantize such a module. In the framework
of the formal quantization scheme such a quantization  is ensured by the construction in the paper \cite{R}.

We do not use this deformation quantization scheme. By contrast,  for certain braided orbits $\Omn$
we  explicitly construct projective modules which play the role  of the cotangent bundles over generic orbits in $gl(m)^*$ in the framework
of the Serre approach. We call these modules {\em cotangent}. Also, we call a quotient $\Omn$ for which this module exists  a
{\em regular} braided variety or a braided {\em  generic} orbit. Specializing our general construction to the case when $R$ is
a super-flip, we get cotangent modules over super-orbits in $gl(m|n)^*$.

In order to construct these cotangent modules we employ the differential calculus on the algebra $\lq$ developed in \cite{GS2}.
In this calculus we do not use any form of the Leibnitz rule which is usually employed in "quantum differential calculus". Instead,
we are dealing with a Koszul type complex (see  section \ref{sec:3}). This approach enables us to compute the
differentials of the functions $\Tr_R L^k,\, k=1,...,m+n$, and to explicitly construct the mentioned
cotangent module provided the eigenvalues $\mu$ and  $\nu$ do not belong to an exceptional set $\E$ which is defined
by zeros of a determinant. So, all quotients $\Omn$ corresponding to the eigenvalues
$(\mu, \nu)\in \C^{\oplus (m+n)}\setminus \E$ are considered to be  braided   generic orbits.

Furthermore, our construction can be extended to appropriate
quotients of the so-called modified REA,  which are in a sense
braided analogs of the enveloping algebras $U(gl(m|n)$). We call
these quotients braided non-commutative (NC) orbits. In section
\ref{sec:5} we present a criterium (similar to that mentioned
above) which ensures  regularity of such an orbit. Considering a
particular case when $R$ is a super-flip, we get a description of
regular or generic NC super-orbits. Note that certain  braided
deformations of generic super-orbits give rise to some Poisson
pencils on these super-orbits which
 will be considered in the subsequent paper \cite{GS3}.

The paper is organized as follows. In the next section we present a short review of the REA and its properties used in the
sequel. In section \ref{sec:3} we present some aspects of the braided differential calculus. It helps us to
formulate a criterium of regularity of algebras we are dealing with. In section \ref{sec:4} we present this criterium in terms
of the quantum eigenvalues. In section \ref{sec:5} we extend these results to the braided NC orbits.
\bigskip

\noindent
{\bf Acknowledgement.} The work of P.S. was partially supported by
the RFBR grant 08-01-00392-a and the joint RFBR and DFG grant 08-01-91953. The work of D.G. and
P.S. was partially supported by the joint RFBR and CNRS grant 09-01-93107.

\section{Reflection Equation Algebra: CH identity and other properties}
\label{sec:2}

Let $V$ be a finite dimensional vector space and $R\in \End(V^{\otimes 2})$ be a Hecke symmetry.
We associate with $R$ two quadratic algebras --- quotients of the free tensor algebra $T(V)$ generated
by the space $V$:
$$
{\rm Sym}_R(V)=T(V)/\langle {\rm Im}(q I-R) \rangle,\quad
{\bigwedge}_R(V)=T(V)/\langle {\rm Im}(\qq I+R) \rangle.
$$
These are $R$-analogs of the usual symmetric and skew-symmetric algebras respectively. Denote
$\Sym^k_R(V)$ and ${\bigwedge}^k_R(V)$ the $k$-th degree homogeneous component of the
algebras in question and introduce the corresponding Hilbert-Poincar\'e series
$$
P_+(t)=\sum_{k\ge 0} t^k \dim\, \Sym^k_R(V),\quad P_-(t)=\sum_{k\ge 0} t^k \dim\, {\bigwedge}^k_R(V).
$$

As was shown in \cite{P}, these series are always rational functions, and, therefore, each of them can be presented as a ratio
of two coprime   polynomials. Denote $m$ (resp., $n$) the degree of the numerator (resp., denominator)
of the rational function $P_-(t)$. Let us call the ordered pair $(m|n)$ {\em bi-rank} of the space $V$ or
the corresponding Hecke symmetry $R$. It is an analog of the super-dimension of the space $V$ and coincides
with it when $R$ is a super-flip or its deformation. In this case $m+n=N=\dim V$. However, in general
it is not so. By using results of  \cite{G} it is possible to construct Hecke symmetries of the bi-rank $(m|n)$
such that $m+n <N$.

In what follows we assume the symmetry $R$ to be skew-invertible. This means that there exists an operator
$\Psi \in \End(\vv)$ such that
\be
\Tr_2 R_{12} \Psi_{23}= \sigma_{13},
\label{def:psi}
\ee
where the (usual) trace is applied to the operator product $R_{12} \Psi_{23}\in \End(V^{\ot 3})$  in the second
space and $\sigma$ is the usual flip. Consider two operators $B:V\to V$ and $C:V\to V$ defined as follows
\be
B=\Tr_1 \Psi,\,\, C=\Tr_2 \Psi.
\label{def:BC}
\ee

These operators play a crucial role in defining $R$-traces mentioned in Introduction. Thus, the operator $C$ comes in the following
way in the formulae
for  the {\em power sums}
$$
 p_k(L)=\Tr_R L^k:=\Tr(L^kC),\quad k\ge 1
$$
(besides, we put $p_0(L):=1$).
As we noticed above these elements belong to the center
of the algebra\footnote{In the sequel we do not need the operator
$B$, it plays an analogous role in constructions related to another version of the REA:
$$
RL_2RL_2 - L_2RL_2R = 0.
$$} $\lq$.

Another family generating  the center $Z(\lq)$ is formed by the so-called Schur functions
(moreover, this family spans $Z(\lq)$ as a vector space).  Any such a
function $\sla(L)$ is associated with a partition $\lambda$ of a non-negative integer  $k$
$$
\la=(\la_1,\, \la_2,...),\quad 0\leq \la_{i+1}\leq \la_i, \quad \sum \la_i =k.
$$
We refer the reader to \cite{GPS1,GPS2} for detailed definition and properties of the Schur functions
for any QMA (see footnote \ref{ftn:4}). Observe that in general the Schur functions are not central, but
span a commutative subalgebra (called {\em characteristic}) of the QMA in question.

It is worth noticing that in any QMA the Schur functions are polynomials in the algebra generators and
satisfy the following multiplication rule
$$
\sla(L) \, s_{\mu}(L)=\sum_{\nu}C_{\la, \mu}^\nu s_{\nu}(L)
$$
where $C_{\la, \mu}^\nu$ are the Littlewood-Richardson coefficients.

The Schur functions corresponding to single-column and single-row partitions $\la=(1^k)$ and $\la=(k)$
respectively, $k=0,1,2,...$, are of special interest. We denote them  $a_k(L)$ and $s_k(L)$ respectively.
The interrelations among the functions from the sets $\{a_k\}$, $\{s_k\}$ and $\{p_k\}$ are described
by the following formulae (below $k\geq 1$)
\begin{eqnarray}
(-1)^k k_q\, a_k(M) + {\textstyle \sum_{r=0}^{k-1}}\, (-q)^{r}a_r(M)\,
p_{k-r}(M) &=& 0,\label{q-anti}\\[1mm]
k_q\, s_k(M) - {\textstyle \sum_{r=0}^{k-1}}\, q^{-r}s_r(M)\,
p_{k-r}(M) &=& 0, \label{q-simm}\\[1mm]
{\textstyle \sum_{r=0}^k}\, (-1)^r a_r(M)\, s_{k-r}(M) &=& 0,
\label{wronski}
\end{eqnarray}
where as usual
$$
k_q:=\frac{q^k-q^{-k}}{q-q^{-1}}.
$$
The relations (\ref{q-anti}) and (\ref{q-simm})  are called  quantum Newton relations.
They differ from the classical versions by factors depending on $q$.
By contrary,  the  relation (\ref{wronski}) (called Wronski one) does not depend on $q$ and  coincides completely
 with  its classical counterpart.

In what follows a distinguished role is played by the following partitions for which we introduce a
special  notation
\begin{eqnarray*}
&&[m|n]:=\bigl( n^m \bigr)\\
&&[m|n]_r := \bigl( n^m, r\bigr)\\
&&[m|n]^k := \bigl( (n+1)^k, n^{m-k}\bigr)\\
&&[m|n]^k_r := \bigl( (n+1)^k, n^{m-k}, r \bigr).
\end{eqnarray*}

In this notation the CH identity for the generating matrix $L$ of the algebra $\lq$ reads (\cite{GPS1})
\be
\sum_{i=0}^{m+n}\Big(\sum_{k=\max(0,i-n)}^{\min(i,m)}(-1)^kq^{2k-i}s_{[m|n]^k_{i-k}}\Big) L^{m+n-i}= 0.
\label{CH0}
\ee
Multiplying this identity by $\schf{[m|n]}(L)$ and taking into account the following quadratic relations among the
Schur functions
\be
\schf{[m|n]}(L)\schf{[m|n]_r^k}(L)=\schf{[m|n]_r}(L)\schf{[m|n]^k}(L)
\label{bil}
\ee
we can rewrite the CH identity (\ref{CH0}) in a factorized form
\be
\Big(\sum_{k=0}^m (-q)^k\,
\schf{[m|n]^k}(L)\, L^{{m-k}}
\Big) \Big(\sum_{r=0}^n q^{-r}\,
\schf{[m|n]_r}(L)\, L^{{n-r}}\Big)=0\, .
\label{CH}
\ee

This  form of the CH identity enables us to introduce the notions of even and odd eigenvalues of the matrix $L$.
Namely, the roots of the first (resp., second) factor in (\ref{CH}) are called {\em even} (resp., {\em odd}) eigenvalues
and are denoted $\mu_i$, $1\leq i \leq m$, (resp., $\nu_i$, $1\leq i \leq n$). Since all coefficients coming in the factors
of the product (\ref{CH}) belong to the center $Z(\lq)$ of the algebra $\lq$, the eigenvalues $\mu_i$ and $\nu_i$ are
treated to be elements of an algebraic extension of the localization $Z(\lq)_{loc}=S^{-1}\,Z(\lq)$ of the
center $Z(\lq)$ by the set $S=\{(\schf{[m|n]}(L))^{k},\, k=1,2,...\}$.

It turns out that all Schur functions can be expressed via the eigenvalues $\mu_i$ and $\nu_i$. Thus, for the
Schur function $\schf{[m|n]}(L)$ we have (see \cite{GPS2})
\be
\schf{[m|n]}(L)\, \mapsto\, s_{[m|n]}(\mu,\nu) = \prod_{i=1}^m\prod_{j=1}^n
\left(q^{-1}\mu_i - q\nu_j\right).
\label{c0}
\ee

In a similar manner we can parameterize the  power sums $p_k(L)=\Tr_R L^k$ (see \cite{GPS4}) :
\begin{equation}
p_k(L)\mapsto p_k(\mu,\nu)=\sum_{i=1}^m d_i \mu_i^k  +
\sum_{j=1}^n  d'_j\nu_j^k\quad \forall\,k\ge 0\,,
\label{pk}
\end{equation}
where  the coefficients $d_i$ and $d'_j$ (called {\em quantum dimensions}) have the form
\be
d_i = q^{-1}\prod_{p=1\atop p\not=i}^m
\frac{\mu_i - q^{-2}\mu_p}{\mu_i-\mu_p}\,
\prod_{j=1}^n \frac{\mu_i - q^2\nu_j}{\mu_i -\nu_j} , \quad
d'_j = \,-q\,\prod_{i=1}^m \frac{\nu_j - q^{-2}\mu_i}{\nu_j-\mu_i}\,
\prod_{p=1\atop p\not=j}^n \frac{\nu_j - q^2\nu_p}{\nu_j-\nu_p}\,.
\label{dd-i}
\ee

{\bf  Example.} Let us consider an example: $m=3, n=2$. In particular, this example covers the case related to
the quantum group $U_q(3|2)$. (As was mentioned in Introduction, in this case  $\dim V = 3+2=5$).

The CH identity (\ref{CH0}) becomes
$$
s_{[3|2]}\, L^5+\left(q^{-1}s_{[3|2]_1} - q s_{[3|2]^1}\right)\,L^4+\left(q^{-2}s_{[3|2]_2} -
s_{[3|2]^1_1} +q^2s_{[3|2]^2}\right)\, L^3+$$
$$
 \left(-q^{-1}s_{[3|2]^1_2} + q s_{[3|2]^2_1} -q^3s_{[3|2]^3}\right)\, L^2 +
\left(s_{[3|2]^2_2} -q^2 s_{[3|2]^3_1}\right)\, L - q\,s_{[3|2]^3_2}\, I= 0\,.
$$

The  bilinear relations (\ref{bil}) read
\begin{eqnarray}
s_{[3|2]}s_{[3|2]^1_1} = s_{[3|2]^1}s_{[3|2]_1}&\hspace*{10mm}&
s_{[3|2]}s_{[3|2]^1_2} = s_{[3|2]^1}s_{[3|2]_2}\nonumber \\
s_{[3|2]}s_{[3|2]^2_1} = s_{[3|2]^2}s_{[3|2]_1}&& s_{[3|2]}s_{[3|2]^2_2} = s_{[3|2]^2}s_{[3|2]_2}\label{bl}\\
s_{[3|2]}s_{[3|2]^3_1} = s_{[3|2]^3}s_{[3|2]_1}&&s_{[3|2]}s_{[3|2]^3_2} = s_{[3|2]^3}s_{[3|2]_2}\,.\nonumber
\end{eqnarray}

Then multiplying the above CH identity by  $s_{[3|2]}$  and employing  these bilinear relations
 we arrive to  the  factorized form of the CH identity:
$$
\left(s_{[3|2]}\,L^3 -q s_{[3|2]^1}\,L^2 +q^2s_{[3|2]^2}\,L-q^3s_{[3|2]^3}\, I\right)
\left(s_{[3|2]}\,L^2 +q^{-1} s_{[3|2]_1}\,L +q^{-2}s_{[3|2]_2}\, I\right)= 0\,.
$$
Thus, we have three even eigenvalues $\mu_1,\,\mu_2,\,\mu_3$ and two odd ones $\nu_1,\,\nu_2$.
In virtue of the Vieta formula we have
$$
q\,\frac{s_{[3|2]^1}}{s_{[3|2]}} =\mu_1+\mu_2+\mu_3,\quad
q^2\,\frac{s_{[3|2]^2}}{s_{[3|2]}} =\mu_1\mu_2+\mu_1\mu_3+\mu_2\mu_3,\quad
q^3\,\frac{s_{[3|2]^3}}{s_{[3|2]}} = \mu_1\mu_2\mu_3,
$$
$$
-q^{-1}\,\frac{s_{[3|2]_1}}{s_{[3|2]}} =\nu_1+\nu_2,\quad
q^{-2}\,\frac{s_{[3|2]_2}}{s_{[3|2]}} =\nu_1\nu_2\,.
$$

So, we can rewrite the CH identity as follows
$$
s_{[3|2]}^2(L-\mu_1)(L-\mu_2)(L-\mu_3)(L-\nu_1)(L-\nu_2)= 0.
$$

Besides, we have the following parametrization of the Schur function  $s_{[3|2]}$:
$$
s_{[3|2]} = \prod_{i=1}^3\prod_{j=1}^2(q^{-1}\mu_i-q\nu_j).
$$

The above parametrization enables us to find a parametrization of all other Schur functions
(see \cite{GPS2} for more detail). For example, we have
$$
s_{[3|2]^3}=q^{-3}\, \mu_1\mu_2 \mu_3  \prod_{i=1}^3\prod_{j=1}^2(q^{-1}\mu_i-q\nu_j),\quad
s_{[3|2]_2}= q^2\, \nu_1\nu_2
\prod_{i=1}^3\prod_{j=1}^2(q^{-1}\mu_i-q\nu_j).
$$

Returning to the general case, we note that in order to realize the localization of the center $Z(\lq)$ by the the set $S$
we should assume that the element $s_{[m|n]}$ does not vanish. In virtue of (\ref{c0}), this requirement entails that
$\mu_i\not=q^2\nu_j$ for all couples $(i,\,j)$. In the next section we shall give a full description of  the
exceptional set $\E$ mentioned in Introduction which includes all families $\mu, \nu$ such that
$\mu_i=q^2\nu_j$ at least for one pair $\mu_i$ and $\nu_j$. Thus, on the complementary set
$\C^{m+n}\backslash \E$ the element $s_{[m|n]}$ is invertible and consequently the localization $Z(\lq)_{loc}=S^{-1}\,Z(\lq)$ is well defined.

\section{Elements of braided differential calculus}
\label{sec:3}

In this section we introduce certain elements of differential calculus on the algebras $\Omn$ based on the
approach suggested in \cite{GPS3} (see  also the references therein). In this approach we use the Koszul type
complexes  whose terms are defined via a series of {\em braided symmetrization} and
{\em skew-symmetrization} projectors. These projectors act in spaces $\L^{\ot k}$ where
$\L={\span}(l_i^j)$ is a vector space generating the algebra $\lq$. Without going into
detail we describe some aspects of this method.

Let $I_- \subset \LL$ be the subspace of $\LL$ spanned by the left hand side of (\ref{RE}). In a sense it is a braided
analog of the usual skew-symmetric subspace. For this reason the algebra $\lq=T(\L)/\langle I_- \rangle$ will
be also denoted  $\Sym_q(\L)$.

Furthermore, in the space $\LL$ there exists another subspace
$$
I_+={\span}(R\,L_1\,R\,L_1+L_1\,R\,L_1\, R^{-1}),
$$
 which  can be considered as an analog of the usual symmetric subspace.
Then the algebra ${\bigwedge}_q(\L)=T(\L)/\langle I_+\rangle$ is an analog of the usual skew-symmetric algebra.
The basic property of the subspaces $I_{\pm}\subset \LL$ is that they are complementary, i.e.
\be
I_+ \bigcap I_-=\{0\},\quad I_+ + I_-=\LL. \label{compl} \ee

Let us suppose that any $k$-th order homogeneous element $f\in \Sym_q^k(\L)=\L^k(R)$ (respectively
$f\in {\bigwedge}_q^k(\L)$) can be presented in the complete "symmetric" (respectively "skew-symmetric") form,
i.e. as an element of the subspace
$$
I_+^{(k)}=I_+\ot \L^{k-2} \bigcap \L\ot I_+ \ot \L^{k-3}\bigcap ... \bigcap \L^{k-2}\ot I_+
$$
$$
 ({\rm respectively},\quad I_-^{(k)}=I_-\ot \L^{k-2} \bigcap \L\ot I_- \ot \L^{k-3}\bigcap ... \bigcap \L^{k-2}\ot I_-).
$$
We call this form {\em canonical}. This presentation of homogeneous elements can be naturally realized via
(skew)symmetrization projectors
$$
P_{\pm}^{(k)}:\L^{\ot k}\to I_{\pm}^{(k)}
$$
with the natural property  $P_{\pm}^{(k)}P_{\pm j} ^{(i)}=P_{\pm}^{(k)}$ where $i+j\leq k+1$ and
 $P_{\pm j}^{(i)}$ stands for  the projector $P_{\pm}^{(i)}$ acting in the product $\L^{\ot k}$  on
the terms  with numbers $j, j+1,...,j+i-1$.

In \cite{GPS3} such projectors have been constructed for $k=2,3$ (also, see formulae (\ref{p+2}) and
(\ref{p+3}) below).  In general, the problem of their explicit construction is still open.

Consider a family of  complexes labelled by positive integers $r$
$$
d: \;{\bigwedge}_q^k(\L)\ot \Sym_q^p(\L) \to {\bigwedge}_q^{k+1}(\L)\ot \Sym_q^{p-1}(\L),
\quad k+p=r
$$
$$
d((y_1 \cdot ...\cdot y_k)\ot (x_1\cdot ...\cdot x_p))=p\, (y_1\cdot ...\cdot y_k\cdot x_1 )
\ot (x_2\cdot ...\cdot x_p).
$$
Here we assume that the elements $y_1\cdot ...\cdot y_k\in {\bigwedge}_q^{k}(\L)$ and
$x_1\cdot ...\cdot x_p\in  \Sym_q^p(\L)=\L^p(R)$ are written in the canonical form. This prevents us
from necessity of using any form of the Leibnitz rule. Note that we have put the factor $p$ in the
above formula for the differential $d$ by analogy with the classical case but in principle it can be replaced by any
non-trivial  factor.

The complexes above can be put together in one complex
\be
d: {\bigwedge}_q^k(\L)\ot \lq\to {\bigwedge}_q^{k+1}(\L)\ot \lq. \label{comp}
\ee
We leave to the reader checking that $d^2=0$. Recall, however, that {\em before} the second
application of the differential $d$ we have to represent the element $y_1\cdot ...\cdot y_k\cdot x_1$
in the canonical form by means of the projector $P_-^{(k+1)}$.

Let us point out that we do not use any transposition between elements from $\lq$ and those
from ${\bigwedge}_q(\L)$. So, we consider the terms of the complex (\ref{comp})
to be one-sided  (namely, right) $\lq$-modules. Note that an attempt to introduce such a
transposition is not compatible with restriction of the space of differential form to braided
orbits (see  \cite{AG} where this problem is discussed on the example of a quantum sphere
(hyperboloid)).

Now, we define the space of first order differentials on the algebras $\Omn$. First, we  apply
the differential $d$ to the elements $\Tr_R L^k$. According to  our scheme, before applying
the differential $d$ to a homogeneous element $\Tr_R L^k$ we have to present it in the canonical
form by means of the symmetrizer $P_+^{(k)}$.

\begin{conjecture}
\label{conj-1}
For the elements   $\Tr_R L^k\in \L^k(R)$ the following relation is valid
\be
P_+^{(k)}\Tr_R L^k= P_{+2}^{(k-1)}  \Tr_R L^k.
\label{tr-symm}
\ee
\end{conjecture}

\begin{proposition} Conjecture \ref{conj-1} is true for the involutive symmetry $R$.
\end{proposition}

Without going into detail we only note that the proof is based on the fact that the element $\Tr_R L^k$
is invariant with respect to the operator $R_{12} R_{23}...R_{k-1\, k}$.

We consider a more interesting and more complicated case when $R$ is not involutive. For the reader's
convenience we reproduce here the explicit formulae for the projectors $P_+^{(2)}$ and $P_+^{(3)}$
from \cite{GPS3}.

We introduce a linear operator ${\sf Q}: {\cal L}^{\otimes 2}\rightarrow {\cal L}^{\otimes 2}$ defining
its action on the basis vectors as follows:
$$
{\sf Q}(L_1\stackrel{.}{\otimes}L_{\bar 2}) = R_{12}L_1\stackrel{.}{\otimes}L_{\bar 2}R_{12}^{-1}.
$$
Here $L_{\bar 2}= R_{12}L_1R_{12}^{-1}$, and the notation $\stackrel{.}{\otimes}$ stands for the
usual matrix product of $L_1$ and $L_{\bar 2}$ but their matrix elements are tensorized in the
resulting matrix instead of being multiplied. The matrix elements of $L_1\stackrel{.}{\otimes}L_{\bar 2}$
form a basis  of the space ${\cal L}^{\otimes 2}$.

Explicitly, this basis is as follows
$$
(L_1\stackrel{.}{\otimes}L_{\bar 2})_{i_1i_2}^{j_1j_2} = l_{i_1}^{a_1}\otimes l_{b_1}^{c_1}R_{a_1i_2}^{b_1a_2}
{R^{-1}}_{c_1a_2}^{j_1j_2}.
$$
(Hereafter, the summation over repeated indices is assumed.) Similarly, as a basis set of the space
${\cal L}^{\otimes 3}$ we shall use the matrix elements of $L_1\stackrel{.}{\otimes}L_{\bar 2}
\stackrel{.}{\otimes}L_{\bar 3}$, with $L_{\bar 3} = R_{23}L_{\bar 2}R^{-1}_{23}$. This basis is
more suitable in working with the REA. For more detail see \cite{GPS3}.  To simplify the writing,
we shall omit the symbol $\stackrel{.}{\otimes}$ in formulae below.

With the above operator ${\sf Q}$ we introduce the projectors $P_\pm^{(2)}$ by the relations
\be
P_+^{(2)} = \frac{1}{2_q^2}\,((q^2+q^{-2})\,{\rm Id} +{\sf Q} + {\sf Q}^{-1}),\qquad
P_-^{(2)} = \frac{1}{2_q^2}\,(2\,{\rm Id} -{\sf Q} - {\sf Q}^{-1}). \label{p+2}
\ee
These operators are complementary in the following sense
$$
P_+^{(2)} + P_-^{(2)} = {\rm Id},\qquad P_\pm^{(2)}P_\pm^{(2)} = P_\pm^{(2)},\qquad
P_\pm^{(2)}P_\mp^{(2)} = 0
$$
(see formula (\ref{compl}) above).
The first relation is evident, the others can be verified with the use of the characteristic polynomial for ${\sf Q}$.

The  projector $P_+^{(3)}:\L^{\ot 3}\to \L^{\ot 3}$ reads
\begin{eqnarray}
P_+^{(3)} &=& \frac{2_q^6}{4\cdot
3_q^2}\,\left(P_{+1}^{(2)}P_{+2}^{(2)}P_{+1}^{(2)}P_{+2}^{(2)}P_{+1}^{(2)}
-a\,P_{+1}^{(2)}P_{+2}^{(2)}P_{+1}^{(2)} +b\,P_{+1}^{(2)}\right)\nonumber\\
&=& \frac{2_q^6}{4\cdot
3_q^2}\,\left(P_{+2}^{(2)}P_{+1}^{(2)}P_{+2}^{(2)}P_{+1}^{(2)}P_{+2}^{(2)}
-a\,P_{+2}^{(2)}P_{+1}^{(2)}P_{+2}^{(2)} +b\,P_{+2}^{(2)}\right)
\label{p+3}
\end{eqnarray}
where
$$
a = (q^4+q^2+4+q^{-2}+q^{-4})/2_q^4,
\quad b = 4_q^2/2_q^8.
$$

\begin{proposition}
For an arbitrary Hecke type symmetry $R$ the formula (\ref{tr-symm}) is valid
for $k=2$ and $k=3$.
\end{proposition}

\noindent{\bf Proof.}
For $k=2$ the claim means that the element $\Tr_R L^2\in \Sym_q^2(\L)$, i.e. it is already symmetrized.
To verify this we calculate the action of $P_+^{(2)}$ on the element $\Tr_RL^2$. For this purpose we use the
following identity which can be easily proved with the use of (\ref{def:psi}) and (\ref{def:BC})
$$
\Tr_R L^2 \equiv \Tr_{R(12)}(L_1L_{\bar 2}R_{12}).
$$
Using the definition of ${\sf Q}^{\pm 1}$ and the cyclic property of the $R$-trace
$$
\Tr_{R(12\dots k)}(R^{\pm 1}_{ii+1}M_{12\dots k}) = \Tr_{R(12\dots k)}(M_{12\dots k}R^{\pm 1}_{ii+1}),
\quad \forall M,\; 1\le i\le k-1,
$$
we get:
$$
 P_+^{(2)}(\Tr_RL^2) =2_q^{-2}\,\Tr_{R(12)}\Big((q^2+q^{-2} +2)L_1L_{\bar 2} R_{12}\Big) = \Tr_RL^2.
$$
Therefore, $\Tr_RL^2\in {\rm Im}(P_+^{(2)})= \Sym_q^2(\L)$.

Turn now to the element $\Tr_RL^3$. Below we use the shorthand notation $R_i:=R_{ii+1}$. Using the
identity
$$
\Tr_RL^3\equiv \Tr_{R(123)}\Big(L_1L_{\bar 2}L_{\bar 3}R_{2}R_{1}\Big)
$$
we first find
\be
P_{+2}^{(2)}(\Tr_RL^3) = 2_q^{-2}\,\Tr_{R(123)}\Big(L_1L_{\bar 2}L_{\bar 3}(2(R_1R_2+R_2R_1)
-\xi R_1+\xi R_1R_2R_1)\Big),
\label{p2-action}
\ee
where $\xi:=q-q^{-1}$.
To obtain the above result we used the cyclic property of $R$-trace and the following
relation
$$
R_iL_{\bar k} = L_{\bar k}R_i, \quad \forall\, i,k: \;i\not=k-1,k.
$$
Now we shall calculate the action $P_+^{(3)}(\Tr_RL^3)$ and compare it with (\ref{p2-action}). As can be
seen from the structure of $P_+^{(3)}$ given in the second line of (\ref{p+3}), we actually have to calculate
the action of $P^{(2)}_{+2}P^{(2)}_{+1}$ and of its square $(P^{(2)}_{+2}P^{(2)}_{+1})^2$ on the right hand
side of (\ref{p2-action}). These actions lead to transformation of $R$-matrix multipliers at $L_1L_{\bar 2}L_{\bar 3}$
under the sign of $R$-trace. We shall not reproduce the calculations in full detail constraining ourselves by
writing down some key intermediate results. First of all we note that the matrix structures
$$
I_A:= R_1R_2R_1 +R_1+R_2,\quad I_B:= R_1R_2+R_2R_1 -\xi (R_1+R_2)
$$
are {\it invariant} with respect to the action of $P_{+1}^{(2)}$ and $P_{+2}^{(2)}$ which precisely means that
$$
P_{+1,2}^{(2)}\Big(\Tr_{R(123)}(L_1L_{\bar 2}L_{\bar 3}\cdot I_{A,B})\Big) =
\Tr_{R(123)}(L_1L_{\bar 2}L_{\bar 3}\cdot I_{A,B}).
$$
In terms of these invariants the transformation of matrix structures under the action of the operator $P_{+1}^{(2)}$
reads:
\begin{eqnarray*}
P_{+1}^{(2)}&:&R_1\rightarrow R_1\\
&& R_2\rightarrow 2_q^{-2}\,(2I_A-\xi I_B - (\xi^2+2)R_1)\\
&&R_1R_2R_1\rightarrow 2_q^{-2}\,((\xi^2+2)I_A +\xi I_B -2 R_1)\\
&&R_1R_2+R_2R_1\rightarrow 2_q^{-2}\,(2\xi I_A +4 I_B +2\xi R_1).
\end{eqnarray*}
The action of $P_{+2}^{(2)}$ are obtained from the above formulae by changing $R_1\leftrightarrow R_2$.

Representing the matrix structure in (\ref{p2-action}) in the form
$$
2(R_1R_2+R_2R_1)
-\xi R_1+\xi R_1R_2R_1\equiv \xi I_A+2I_B+\xi R_2
$$
we find
$$
 P_{+2}^{(2)}P_{+1}^{(2)}:\quad\xi I_A+2I_B+\xi R_2\;\rightarrow\;
\frac{2_q^4+4}{2_q^4}\,\xi I_A+\frac{(2_q^4-\xi^2)}{2_q^4}\,2I_B
+\frac{(2_q^2-2)^2}{2_q^4}\,\xi R_2.
$$
The action of $(P_{+2}^{(2)}P_{+1}^{(2)})^2$ leads to more cumbersome expression:
\begin{eqnarray*}
 (P_{+2}^{(2)}P_{+1}^{(2)})^2:&&\xi I_A+2I_B+\xi R_2\;\rightarrow\\
&&\hspace*{-15mm} \Big(1+\frac{8(2+2_q^2(2_q^2-2))}{2_q^8}\Big)\,\xi I_A
+\Big(1-\frac{2\xi^2 (2+2_q^2(2_q^2-2))}{2_q^8}\Big)\,2I_B +\frac{(2_q^2-2)^4}{2_q^8}\,
\xi R_2.
\end{eqnarray*}
To get the action of the operator $P_+^{(3)}$ we should compose the linear combination in accordance
with (\ref{p+3}):
$$
P_+^{(3)}(\Tr_RL^3) = \frac{2_q^6}{4\cdot 3_q^2}\,((P_{+2}^{(2)}P_{+1}^{(2)})^2 -
a\,P_{+2}^{(2)}P_{+1}^{(2)} +b)\,P_{+2}^{(2)}(\Tr_RL^3).
$$
On taking into account (\ref{p2-action}), the above results for the action of powers of $P_{+2}^{(2)}P_{+1}^{(2)}$
together with the values of coefficients $a$ and $b$, we come to the formula (\ref{tr-symm}).
\hfill\rule{6.5pt}{6.5pt}
\medskip

Note that in order to prove the conjecture for  $k\geq 4$ we need the explicit  form of the projectors $P_+^{(k)}$.

As follows from the conjecture above the result of applying $d$ to the element $\Tr_R L^k$ equals
$k\,P_{+2}^{(k-1)}  d_1(\Tr_R L^k)$ where $d_1$ stands for the differential (\ref{comp}) applied to the
first factor. Whereas the other factors are assumed to be symmetrized via the projector $P^{(k-1)}_{+2}$.
However, this projector commutes with $d_1$. Therefore, we can apply the operator $d_1$ first and then
apply the symmetrizer $P_{+2}^{(k-1)}$ to the result.

Upon writing $\Tr_RL^k$ in the explicit form
$$
\Tr_R L^k=l_{p_1}^{p_2}l_{p_2}^{p_3}... l_{p_{k-1}}^{p_{k}}C^{p_{1}}_{p_{k}}
$$
we have
$$
d_1(\Tr_R L^k)=d(l_{p_1}^{p_2})\ot  l_{p_2}^{p_3}... l_{p_{k-1}}^{p_{k}}C^{p_{1}}_{p_{k}}.
$$
(Here all indices run from 1 till $N=\dim V$.) Though in this writing we do not apply the  projector $P_{+2}^{(k-1)}$,  it does not affect
the element $d_1(\Tr_R L^k)$  but only its presentation.

Now, we are able to introduce the space  $\Om^1(\O_{\mu,\nu})$ of differential 1-forms on the algebra $\Omn$
as the quotient
$$
{\Om^1(\O_{\mu,\nu})}={\bigwedge}_q^1(\L)\ot \Omn /\langle d(\Tr_R L^1),...,d(\Tr_R L^{m+n}) \rangle
$$
where the denominator is treated to be a submodule of the right $\Omn$-module ${\bigwedge}_q^1(\L)\ot \Omn$.

Though we do not need spaces of higher order differential forms, we want to mention their construction.
Thus, we define the space $\Om^k(\O_{\mu,\nu})$ to be the quotient of the right $\Omn$-module
${\bigwedge}_q^k(\L)\ot \Omn$ over submodule generated by elements
$\omega\cdot d(\Tr_R L^i), \, i=1,...,m+n$ where $\omega$ runs over the space ${\bigwedge}_q^{k-1}(\L)$.
Note that the first $k$ factors of the elements $\omega\cdot d(\Tr_R L^i)$
must be presented in the canonical form, i.e.  skew-symmetrized.

\section{Cotangent modules over braided orbits}
\label{sec:4}
In the classical case  all  spaces $\Om^k(\O_{\mu,\nu})$ constructed in the previous section are vector bundles
provided that the corresponding orbit is generic. In what follows we restrict ourselves to the vector bundle
of 1-forms and describe (the space of sections of) this bundle in the spirit of the Serre approach as a projective module over the
corresponding coordinate ring. We shall call it the {\it cotangent} module. More precisely, we construct such a
module  over the algebra $\Omn$ for generic $\mu$ and $\nu$. The values of these
parameters for which the cotangent module does not exist will be included in an exceptional set $\E$.

Let us set
$$
a^i_j(1)=C_j^i,\qquad a^i_j(k)= l_j^{p_1}\, l_{p_1}^{p_2}... l_{p_{k-2}}^{p_{k-1}}C^i_{p_{k-1}},\quad k\geq 2,\,\, 1\leq i,j \leq N.
$$
Then we have
$$
d(\Tr_R L^k)=(d\,l_i^j)\ot a^i_j(k).
$$

Consider a matrix
$$
A=\big(a(1),\, a(2),..., a(m+n)\big)
$$
where $a(k),\,1\leq k \leq m+n $ is the column
$$
\big(a_1^1(k), a_1^2(k),...,a_1^N(k),a_2^1(k), a_2^2(k),...,a_2^N(k),..., a_N^1,a_N^2(k),...,a_N^N\big)^t.
$$
Hereafter $t$ stands for the transposition. The column $a(k)$ is treated to be the "gradient" of the power
sum $\Tr_R L^k$. Thus, the size of the matrix $A$ is $N^2\times (m+n)$.

Consider another matrix $B$ of the size $(m+n)\times N^2$ defined as follows
$$
B=\big(b(1), b(2),...,b(m+n)\big)^t
$$
where the row $b(k)$ reads
$$
b(k)=\big(b_1^1(k), b^1_2(k),...,b^1_N(k),b^2_1(k), b_2^2(k),...,b^2_N(k),..., b^N_1,b^N_2(k),...,b_N^N\big)
$$
and
$$
b_i^j(1)=\delta_i^j, \quad b_i^j(k)=l_i^{p_1}\, l_{p_1}^{p_2}... l_{p_{k-2}}^j,\quad k\geq 2,\,\, 1\leq i,j \leq N
$$

Now, we calculate the matrix product $B\cdot A$ :
 $$
B\cdot A =
\left(\begin{array}{cccc}
\Tr_R I&\Tr_R L&...&\Tr_R L^{m+n-1}\\
\Tr_R L&\Tr_R L^2&...&\Tr_R L^{m+n}\\
\Tr_R L^2&\Tr_R L^3&...&\Tr_R L^{m+n+1}\\
...&...&...&...\\
\Tr_R L^{m+n-1}&\Tr_R L^{m+n}&...&\Tr_R L^{2(m+n-1)}
\end{array}\right).
$$

Formula (\ref{sys1}) allows us to express the entries of this matrix via the eigenvalues $\mu$ and $\nu$
\be
B\cdot A =
\left(\begin{array}{cccc}
p_0(\mu,\nu)&p_1(\mu,\nu)&...&p_{m+n-1}(\mu,\nu)\\
p_1(\mu,\nu)&p_2(\mu,\nu)&...&p_{m+n}(\mu,\nu)\\
p_2(\mu,\nu)&p_3(\mu,\nu)&...&p_{m+n+1}(\mu,\nu)\\
...&...&...&...\\
p_{m+n-1}(\mu,\nu)&p_{m+n}(\mu,\nu)&...&p_{2(m+n-1)}(\mu,\nu)
\end{array}\right).
\label{mat}
\ee

Now, we are going to calculate the determinant of  the matrix (\ref{mat}) and to find the conditions
under which this matrix is invertible. Using relations (\ref{pk}) we can factorize this matrix into the
product of two square $(m+n)\times (m+n)$ matrices
$$
\left(\begin{array}{cccccc}
d_1&...&d_m&d'_1&...&d'_n\\
d_1\mu_1&...&d_m\mu_m&d'_1\nu_1&...&d'_n\nu_n\\
d_1\mu_1^2&...&d_m\mu_m^2&d'_1\nu_1^2&...&d'_n\nu_n^2\\
...&...&...&...&...\\
d_1\mu_1^{m+n-2}&...&d_m\mu_m^{m+n-2}&d'_1\nu_1^{m+n-2}&...&d'_n\nu_n^{m+n-2}\\
d_1\mu_1^{m+n-1}&...&d_m\mu_m^{m+n-1}&d'_1\nu_1^{m+n-1}&...&d'_n\nu_n^{m+n-1}
\end{array}\right)
\left(\begin{array}{ccccc}
1&\mu_1&\mu_1^2&...&\mu_1^{m+n-1}\\
...&...&...&...&...\\
1&\mu_m&\mu_m^2&...&\mu_m^{m+n-1}\\
1&\nu_1&\nu_1^2&...&\nu_1^{m+n-1}\\
...&...&...&...&...\\
1&\nu_n&\nu_n^2&...&\nu_n^{m+n-1}\\
\end{array}\right)\,.
$$
Now, on taking out the factors $d_i$ and $d'_i$ from the determinant of the first matrix we get that the
determinant of the matrix (\ref{mat}) equals
$$
\det (B\cdot A) = \prod_{i=1}^m d_i\, \prod_{j=1}^n d'_j \,\Big(\prod_{i<j}(\mu_i-\mu_j)\,
\prod_{i,j}(\mu_i-\nu_j) \,\prod_{i<j}(\nu_i-\nu_j)\Big)^2.
$$
Here we used the formula for the determinant of a Wronski matrix.

Again, by means of formulae  (\ref{dd-i}), we conclude that the matrix (\ref{mat}) is invertible iff
\be
\mu_i\not=q^2\mu_j,\quad \nu_i\not=q^2\nu_j,\quad \mu_i\not=q^2\nu_j\, \quad 1\leq i \leq m,\, 1\leq j \leq n.
\label{cond}
\ee
Let us unite all values of  $(\mu, \nu)$ which do not satisfy this condition into the set $\E$.

Now, for all values of the parameters which do not belong to $\E$ we can construct the cotangent module
on the algebra $\Omn$. Let us denote $C = (B\cdot A)^{-1}$ and consider the $N^2\times N^2$ matrix
$\oe=A\cdot C\cdot B$. It can be easily seen that the matrix $\oe$ is an idempotent: $\oe^2=\oe$. Thus, the
right $\Omn$-module
$$
\oM=\oe\,\Omn^{\oplus N^2}
$$
 is projective. It is generated by the "gradients" of the power sums $\Tr_R L^k, 1\leq k\leq m+n$. The complementary
module
$$
M=e\,\Omn^{\oplus N^2},\quad {\rm where} \quad  e=I-\oe
$$
gives an explicit realization of the space $\Om^1(\O_{\mu,\nu})$ as a projective module.

Concluding this section, we want to make the following observation.
Any generic orbit in $gl(m)^*$ can be given via different (but equivalent) systems of equations. Thus, instead of the system
 (\ref{sys}) we can consider that obtained by fixing the values of
the coefficients of the CH identity (namely, the elementary symmetric polynomials in eigenvalues).
In a general case we have a similar situation:  fixing  values of the eigenvalues of the generating matrix $L$ is equivalent to using
the system (\ref{sys1}).

In order to show this equivalence we first express  all coefficients of the  CH polynomial  (the leading coefficient included)
via the power sums $\Tr_R L^k, \, 1\leq k \leq m+n$.  To this end we express the elementary symmetric functions $a_r(L),\, 0\leq r \leq m+n$
via these sums by employing  formulae (\ref{q-anti}). Then by using the Jacobi-Trudy formulae (see
 \cite{FH}) we can express the coefficients $[m|n]^k$ and $[m|n]_r$ of the factorized CH identity  (\ref{CH}) via the functions
$a_r(L),\, 0\leq r \leq m+n$. Thus, all coefficients of the CH identity (\ref{CH0}) can be realized as some polynomial
expressions in the power sums in question. It remains to note that any elementary symmetric (in the usual sense)
polynomial in the roots of the CH identity can be presented as the corresponding  coefficient of the CH polynomial
divided by $s_{[m|n]}$ and consequently, as a rational function in the power sums $\Tr_R L^k$, $1\leq k \leq m+n$.

Nevertheless, namely the system (\ref{sys1}) is the most convenient for constructing the cotangent projective modules over braided generic orbits.

\section{Extension to  braided NC orbits}
\label{sec:5}

Besides the REA $\lq$  there are known other quantum matrix algebras with similar
properties of generating matrix $L$. For these algebras certain quotients looking like
braided orbits can also be defined. The well known example is the algebra $U(gl(m|n))$. Its generating
matrix satisfies a CH identity with central coefficients.  Also, a formula analogous to (\ref{pk}) is valid.
Thus, technique developed in the previous section can be applied for definition of analogs of generic
orbits in $U(gl(m|n))$.

The simplest way to realize this program is to pass to the so-called modified Reflection Equation
Algebra (mREA).
The defining relations of mREA $\lhq$ are similar to that of $\lq$ (\ref{RE}) but with linear terms in the
right hand side:
\be
R\,\hL_1\,R\,\hL_1-\hL_1\,R\,\hL_1\, R=\h(R\,\hL_1-\hL_1\,R),
\label{mRE}
\ee
where $\hL=\|\hl_i^{j}\|$, $1\leq i,\,j\leq \dim\,V$ and $\hL_1=\hL\ot 1$.
All objects related to the mREA $\lhq$ will be denoted by hatted letters.

We introduced a parameter $\h$ in the definition of the mREA in order to present this algebra as a deformation of
$\lq$ in the case when the Hecke symmetry $R=R(q)$ is a deformation of the super-flip
 $R(1)\in\End(\vv)$ where $V$ is a super-space of super-dimension $(m|n)$. In this case
the mREA turns into the algebra $U(gl(m|n)_\h)$ as $q\to 1$
(the subscript $\h$ means that we have introduced the factor
$\h$ in the Lie bracket of the super-Lie algebra $gl(m|n)$).
 For this reason we treat the algebra
$\lqh$ to be a braided analog of the enveloping algebra $U(gl(m|n)_\h)$.

Observe that for $q\not=1$ the algebras $\lq$ and $\lqh$ are isomorphic to each other
(though it is not so for the algebras $\Sym(gl(m|n))$ and $U(gl(m|n)_\h)$). In order to construct their isomorphism we put
\be
L=\hL-{{\h}\over{\xi}},\quad {\rm where}\quad \xi=q-\qq
\label{shift}
\ee
Then,  the system (\ref{RE}) turns into that (\ref{mRE}).
However, this isomorphism fails as $q\to 1$.

Now, we state that the matrix $\hL$ obeys the CH identity
\be
 \sum_{k=1}^{m+n} \hc_k(\hL) \hL^k=0,
\label{odin}
\ee
with central coefficients: $\hc_k\in Z(\lqh)$. In order to find the corresponding  CH polynomial we should make the
shift (\ref{shift}) in the CH (\ref{CH-ident}) and reduce the resulting expression
$$
\sum_{k=1}^{m+n} c_k(\hL-{{\h}\over{\xi}}) (\hL-{{\h}\over{\xi}})^k=0
$$
to the form (\ref{odin}).

By straightforward but tedious computations  it is possible to show that
 the  coefficients $\hc_k(\hL)$ of the  polynomial in (\ref{odin}) have a finite limit as
$q\to 1$. (Note that in the case $n=0$ this property was proven in \cite{GS1}.) Thus, by passing to the limit $q\to 1$ we
get the CH identity for the matrix $\hL$  generating the algebra  $U(gl(m|n)_\h)$
such that the coefficients of this identity are central polynomials  in the generators of the algebra in question.

Denote $\hm_i,\,1\leq i\leq m$, and $\hn_j,\, 1\leq j\leq n$, the roots of the equation
 \be
\sum_{k=1}^{m+n} \hc_k(\hL) t^k=0 \label{NCCH}
\ee
corresponding respectively to $\mu_k$ and $\nu_k$. Namely, we have $\hm_k=\mu_k+{{\h}\over{\xi}}$,
$\hn_k=\nu_k+{{\h}\over{\xi}}$. The roots $\hm_k$ and $\hn_k$ are called respectively {\em even} and {\em odd }
eigenvalues of the matrix $\hL$.

Expressing the power sums $\hp_k(\hL)=\Tr_R \hL^k$ via these eigenvalues we get the  formula analogous to
(\ref{pk}) but with different expressions for quantum dimensions:
\be
\hd_i = q^{-1}\prod_{p=1\atop p\not=i}^m
\frac{\hm_i - q^{-2}\hm_p-\qq\h}{\hm_i-\hm_p}\,
\prod_{j=1}^n \frac{\hm_i - q^2\hn_j+q\h}{\hm_i -\hn_j}\, , \label{newd} \ee
\be
\hd'_j = -\,q\,\prod_{i=1}^m \frac{\hn_j - q^{-2}\hm_i-\qq\h}{\hn_j-\hm_i}\,
\prod_{p=1\atop p\not=j}^n \frac{\hn_j - q^2\hn_p+q\h}{\hn_j-\hn_p}\,. \label{newdd} \ee

In order to prove these formulae it suffices to observe that
$$
\Tr_R f(L)=\sum_{i=1}^m f(\mu_i) d_i+\sum_{j=1}^n f(\nu_i)d'_i
$$
where $d_i$ and $d'_j$ are defined by (\ref{dd-i}) and $f(t)$ is an arbitrary polynomial.

Taking the limit $q\to 1$ in the CH (\ref{NCCH}), we get a formula
for the power sums in the algebra $U(gl(m|n))$. Namely, we obtain that in this algebra the quantum
dimensions are
$$
\hd_i = \prod_{p=1\atop p\not=i}^m
\frac{\hm_i - \hm_p-\h}{\hm_i-\hm_p}\,
\prod_{j=1}^n \frac{\hm_i - \hn_j+\h}{\hm_i -\hn_j}\, ,$$
$$
\hd'_j = -\,\,\prod_{i=1}^m \frac{\hn_j - \hm_i-\h}{\hn_j-\hm_i}\,
\prod_{p=1\atop p\not=j}^n \frac{\hn_j - \hn_p+\h}{\hn_j-\hn_p}\,.
$$

Going back to the general case we consider the following quotients of the algebras $\lqh$
$$
\hOmn=\lqh/\langle \Tr_R L-\hp_1(\mu, \nu),..., \Tr_R L^{m+n}-\hp_{m+n}(\mu, \nu)\rangle ,
$$
where the functions $\hp_k(\mu, \nu)$ are defined by (\ref{pk}) but with quantum dimensions given by
(\ref{newd}) and (\ref{newdd}). These quotient are called {braided NC orbits}.

Let us define the projective module $\hMmn$ similar to $\Mmn$. We set $\hMmn=\he \hOmn^{\oplus N^2}$
 where $\he$ is defined by a formula similar to that for $e$. The only modification consists in
defining the exceptional set $\hE$ of the values $\hm, \hn$ for which  construction of the module $\hMmn$ fails.
The set $\hE$ contains all parameters $\hm, \hn$ for which at least one of the following conditions fails is not fulfilled
$$
\hm_i - q^{-2}\hm_j-\qq\h\not=0,\quad \hn_j - q^2\hn_j+q\h\not=0,\quad\hm_i - q^2\hn_j+q\h\not=0.
$$

By analogy with the previous case, the  $\hOmn$-module $\hMmn$ is called {\em cotangent} one. Upon taking
the limit $q\to 1$ we get the {\em cotangent module} over a {\em NC super-orbit}. The corresponding exceptional set is a
specialization  of $\hE$ where we put $q=1$.

In conclusion we would like to emphasize that the family of regular orbits in a classical (or super-) case
 is bigger than in the case of a braided deformation.
For instance, compare this family for the classical  case $gl(2)$ and that for its braided (NC) counterpart. If in the former case
the only restriction
on the eigenvalues is $\mu_1\not=\mu_2$ in the latter case there are two restrictions
$\hm_1\not=q^2 \hm_2+ \qq \h$ and $\hm_2\not=q^2 \hm_1+ \qq \h$.
In general, they coincide with each other iff $q=1$ and $\h=0$.

\end{document}